\documentclass[12pt]{article}

\usepackage{amsmath}
\usepackage{amssymb}

\newcommand\E{\ensuremath{\mathbb{E}}}

\newcommand{\bone}{\boldsymbol{1}}

\newcommand{\tzero}{\textit{0}}

\title{A note on the state occupancy distribution for Markov chains}
\author{Phil Pollett}

\begin{document}
\maketitle

\begin{center}
\begin{minipage}{14cm}
{\bf Abstract} \ \ 
In a recent paper, Shah~\cite{Sha24} derived an explicit expression for
the distribution of occupancy times for a two-state Markov chain, using
a method based on enumerating sample paths. We consider here the more
general problem of finding the distribution of occupancy times for
countable-state Markov chains in discrete time. Our approach, which
employs generating functions, leads to arguably simpler formulae for 
the occupancy distribution for the two-state chain.
\end{minipage}
\end{center}


\section{The general setup}
Let $(X_n,\, n\geq 0)$ be a discrete-time Markov chain taking values in
a countable state space $S$ with 1-step transition matrix $P=(p_{ij},\,
i,j\in S)$. Let $U$ be some subset of $S$. We are interested here in
$N_n$, the number of times the chain visits $U$ up to time $n$:
\begin{equation}
N_n = \sum_{m=1}^n 1_{\{X_m\in U\}},\qquad n\geq 1.
\label{Nn}
\end{equation}
Our approach is similar to that outlined in Section~3.1 of~\cite{Ser00}. 
However, note that here the initial state $X_0$
does not contribute to the sum.

\medskip
Let $g_i(n,k)=\Pr(N_n=k|X_0=i)$, for $n\geq 1$ and $0\leq k\leq n$. 
First observe that
$g_i(1,1) =\sum_{j\in U} p_{ij}$,
and $g_i(1,0) = 1- g_i(1,1) =\sum_{j\in U^c} p_{ij}$.
By conditioning of the first state visited, we get, for $n\geq 2$,
$$
g_i(n,k)
=\sum_{j\in U} p_{ij} g_j(n-1,k-1)
+\sum_{j\in U^c} p_{ij} g_j(n-1,k),
\qquad 1\leq k\leq n-1,
$$
and
$$
g_i(n,n) =\sum_{j\in U} p_{ij} g_j(n-1,n-1),
\qquad
g_i(n,0) =\sum_{j\in U^c} p_{ij} g_j(n-1,0).
$$
It will be convenient to use matrix notation, writing 
$$
P = \left(
\begin{matrix}
P_{UU} &P_{UU^c} \\
P_{U^cU} &P_{U^cU^c}
\end{matrix}
\right)
$$
and $g(n,k)=(g_i(n,k),\, i\in S)$, a column vector, as
$$
g(n,k) = \left(
\begin{matrix}
g_{U}(n,k)  \\g_{U^c}(n,k)
\end{matrix}
\right),
$$
and then setting
$$
A=\left(
\begin{matrix}
P_{UU} &\tzero\\
P_{U^cU} &\tzero
\end{matrix}
\right),\quad
B= \left(
\begin{matrix}
\tzero &P_{UU^c}\\
\tzero &P_{U^cU^c}
\end{matrix}
\right).
$$
Then, for $n\geq 2$,
\begin{equation}
g(n,k) = A g(n-1,k-1) + B g(n-1,k),
\qquad 1\leq k\leq n-1,
\label{gnk}
\end{equation}
and, with $g(1,0) = B\bone$ and $g(1,1) = A\bone$, where $\bone$
is a column vector of $1$s,
\begin{equation}
g(n,n) = A g(n-1,n-1)=A^n\bone, \quad g(n,0) = B g(n-1,0)=B^n\bone.
\label{gn0}
\end{equation}
As we shall see, this form lends itself to further analysis, taking
generating functions on the first entry, but it is also convenient for
numerical calculations, especially if using a package that exploits
matrix computations; Matlab was used to verify the formulae given
in Theorem~1 below.

Now define the generating function $G_i(\cdot, k)$, for $k\geq 0$, by
$$
G_i(t,k)=\sum_{n=k}^\infty g_i(n,k) t^{n-k},
\qquad |t|<1,
$$
where for convenience $g_i(0,0)=1$, 
and $G(t,k)$ the column vector $(G_i(t,k),\, i\in S)$:
$$
G(t,k)=\sum_{n=k}^\infty g(n,k) t^{n-k},
\qquad |t|<1.
$$
For example, since $|t|<1$,
$G(t,0) =\sum_{n=0}^\infty B^n t^n \bone =(I-Bt)^{-1}\bone$.
Then, using (\ref{gnk}) and (\ref{gn0}), it is not difficult to see that
$$
G(t,k) = AG(t,k-1) + BtG(t,k),\qquad  k\geq 1.
$$
Therefore, since $|t|<1$,
\begin{equation}
G(t,k)=(I-Bt)^{-1} AG(t,k-1)=[(I-Bt)^{-1} A]^k (I-Bt)^{-1}\bone,
\label{Gtk}
\end{equation}
this being valid also for $k=0$.

We can proceed a little further using the identity $(I-Bt)^{-1}
=\sum_{n=0}^\infty B^n t^n$. Since
$$
B^n=
\left(
\begin{matrix}
\tzero &P_{UU^c} P_{U^cU^c}^{n-1}\\
\tzero &P_{U^cU^c}^n
\end{matrix}
\right),
\qquad
$$
the ingredients of (\ref{Gtk}) may be written out as
\begin{equation}
(I-Bt)^{-1}
=
\left(
\begin{matrix}
I &tP_{UU^c} (I-tP_{U^cU^c})^{-1}\\
\tzero &(I-tP_{U^cU^c})^{-1}
\end{matrix}
\right),
\label{IBt}
\end{equation}
and
$$
[(I-Bt)^{-1}A]^k
=
\left(
\begin{matrix}
V &\tzero \\
W &\tzero
\end{matrix}
\right) V^{k-1},
$$
where $V=P_{UU}+tP_{UU^c}(I-tP_{U^cU^c})^{-1}P_{U^cU}$, and
$W=(I-tP_{U^cU^c})^{-1}P_{U^cU}$.

\section{The two-state chain}

This has $S=\{0,1\}$ with
$$
P= \left(
\begin{matrix}
1-p &p\\
q   &1-q
\end{matrix}
\right),
$$
where $p,q\in (0,1)$. 
The following result specifies the occupancy distribution for State~0 
(the occupancy distribution for State~1 can then be obtained simply
by relabling the states).

\medskip
\noindent
{\bf Theorem}~1. 
If $p+q=1$, then, for both $i$,
$$
g_i(n,k) = \binom{n}{k} q^k p^{n-k},
\qquad 0\leq k \leq n.
$$
Otherwise, set $r=1-p-q$. Then,
(a) $g_0(n,n)=(1-p)^n$, and, for $n\geq 1$ and $0\leq k \leq n-1$,
$$
g_0(n,k)
=p(1-q)^{n-(2k+1)}\!\!\!\!\sum_{j=0}^{\min\{k,\,n-k\}}
\binom{n-k}{j}\binom{n-j}{k-j} 
\left(1-\frac{j}{n-k} q\right)(1-q)^{k-j}(1-p)^{k-j}(-r)^j.
$$
(b) $g_1(n,0)=(1-q)^n$, and, for $n\geq 1$ and $1\leq k\leq n$,  
$$
g_1(n,k)=q(1-p)^{2k-1-n}\!\!\!\!\sum_{j=0}^{\min\{k,\,n-k\}}
\binom{k}{j}\binom{n-j}{n-k-j} 
\left(1-\frac{j}{k} p\right)(1-p)^{n-k-j}(1-q)^{n-k-j}(-r)^j.
$$

\noindent
{\em Proof\/}. Let $U=\{0\}$. Then,
$P_{UU}=1-p$, $P_{UU^c}=p$, $P_{U^cU}=q$, $P_{U^cU^c}=1-q$,
and hence
$$
A=\left(
\begin{matrix}
1-p &0\\
q   &0
\end{matrix}
\right),
\quad
B= \left(
\begin{matrix}
0 &p\\
0 &1-q
\end{matrix}
\right).
$$
From (\ref{IBt}) we get
$$
(I-Bt)^{-1}
=
\frac{1}{1-(1-q)t}
\left(
\begin{matrix}
1-(1-q)t &pt\\
0 &1
\end{matrix}
\right),
$$
and so
\begin{equation}
G(t,0)
=
\frac{1}{1-(1-q)t}
\left(
\begin{matrix}
1-rt \\
1
\end{matrix}
\right).
\label{Gt01}
\end{equation}
Furthermore,
$$
(I-Bt)^{-1}A 
= 
\frac{1}{1-(1-q)t}
\left(
\begin{matrix}
1-p-rt  &0\\
q       &0
\end{matrix}
\right),
$$
and so, from (\ref{Gtk}),
\begin{equation}
G(t,k)
= 
\frac{(1-rt)(1-p-rt)^{k-1}}{(1-(1-q)t)^{k+1}}
\left(
\begin{matrix}
1-p-rt\\
q     
\end{matrix}
\right),
\qquad k\geq 1.
\label{Gtk1}
\end{equation}
We may now invert (\ref{Gt01}) and (\ref{Gtk1}) using the 
binomial series: for $k\geq 1$ and $|a|<1$,
$$
\frac{1}{(1-a)^k}=\sum_{i=0}^\infty \binom{k+i-1}{i} a^i.
$$
First, if $p+q=1$, that is $r=0$, then
$$
G(t,k)
= 
\frac{q^k}{(1-pt)^{k+1}}
\left(
\begin{matrix}
1\\
1     
\end{matrix}
\right),
\qquad k\geq 0.
$$
It follows that, for both $i$,
$$
g_i(n+k,k) = \binom{k+n}{n} q^k p^n,
\qquad k\geq 0,
$$
and hence that, for $n\geq 0$ and $k=0,1,\dots, n$,
$$
g_i(n,k) = \binom{n}{k} q^k p^{n-k}.
$$

Now suppose that $r\neq 1$.  From (\ref{Gt01}) we get
$g_0(0,0)=1$, $g_0(n,0)=p(1-q)^{n-1}$, $n\geq 1$, and
$g_1(n,0)=(1-q)^n$, $n\geq 0$.
To deal with the remaining cases, let us focus on $G_1$,
first observing that
$$
\frac{q}{(1-(1-q)t)^{k+1}} = \sum_{i=0}^\infty b_i(k) t^i,
\quad \text{where} \quad b_i(k)=\binom{k+i}{i} q(1-q)^i,
$$
and
\begin{multline*}
(1-rt)(1-p-rt)^{k-1}
=
(1-rt)\sum_{i=0}^{k-1} \binom{k-1}{i} (-r)^i (1-p)^{k-1-i} t^i\\
=
\sum_{i=0}^{k-1} \binom{k-1}{i} (-r)^i (1-p)^{k-1-i} t^i
+\sum_{i=1}^{k} \binom{k-1}{i-1} (-r)^i (1-p)^{k-i} t^i\\
=
(1-p)^{k-1}+ 
\sum_{i=1}^{k-1} \left( \binom{k-1}{i} 
+\binom{k-1}{i-1} (1-p)\right) (1-p)^{k-i-1} (-r)^i t^i
+ (-r)^k t^k\\
=
(1-p)^{k-1}+ 
\sum_{i=1}^{k-1} \left( \binom{k}{i} 
-p\binom{k-1}{i-1}\right) (1-p)^{k-i-1} (-r)^i t^i
+ (-r)^k t^k.
\end{multline*}
Let us write
$$
(1-rt)(1-p-rt)^{k-1} = \sum_{i=0}^k a_i(k) t^i,
$$
where $a_0(k)=(1-p)^{k-1}$, $a_k(k)=(-r)^k$, and, for $1\leq i\leq k-1$,
$$
a_i(k) 
= \left( \binom{k}{i} -p\binom{k-1}{i-1}\right) (1-p)^{k-1-i} (-r)^i
= \binom{k}{i}\left(1-\frac{i}{k} p\right)(1-p)^{k-1-i} (-r)^i,
$$
noting that the latter is now valid also for $i=0$ and $i=k$.
Then, for $k\geq 1$,
$G_1(t,k) = \sum_{i=0}^\infty c_i(k) t^i$, where
$c_i(k)=\sum_{j=0}^{\min\{i,\,k\}} a_j(k) b_{i-j}(k)$. Therefore,
\begin{multline*}
c_i(k)
=
\sum_{j=0}^{\min\{i,\,k\}} 
\binom{k}{j}\left(1-\frac{j}{k} p\right)(1-p)^{k-1-j} (-r)^j
\binom{k+i-j}{i-j} q(1-q)^{i-j}\\
=q(1-p)^{k-1-i}\!\!\!\sum_{j=0}^{\min\{i,\,k\}}
\binom{k}{j}\binom{k+i-j}{i-j} 
\left(1-\frac{j}{k} p\right)(1-p)^{i-j}(1-q)^{i-j}(-r)^j.
\end{multline*}
Of course $c_i(k)=g_1(i+k,k)$, $k\geq 0$, $i \geq 0$, and so
$g_1(n,k)=c_{n-k}(k)$, $k\geq 0$, $n\geq k$.
We deduce that, for $n\geq 1$ and $1\leq k\leq n$,  
$$
g_1(n,k)=q(1-p)^{2k-1-n}\!\!\!\!\sum_{j=0}^{\min\{k,\,n-k\}}
\binom{k}{j}\binom{n-j}{n-k-j} 
\left(1-\frac{j}{k} p\right)(1-p)^{n-k-j}(1-q)^{n-k-j}(-r)^j.
$$
This finishes the proof of Part~(a).

To obtain $g_0(n,k)$ we could use the fact that
$G_0(t,k)=(1-p-rt)q^{-1} G_1(t,k)$.
Alternatively, noting that the occupancy of State~1 is $n-N_n$,
we may simply swap $p$ and $q$ and replace $k$
by $n-k$ in our formulae for $g_1(n,k)$.
For example, $g_0(n,n)=(1-p)^n$ since $g_1(n,0)=(1-q)^n$.
More generally, for $n\geq 1$ and $k=0,1,\dots, n-1$,
$$
g_0(n,k)
=p(1-q)^{n-(2k+1)}\!\!\!\!\sum_{j=0}^{\min\{k,\,n-k\}}
\binom{n-k}{j}\binom{n-j}{k-j} 
\left(1-\frac{j}{n-k} q\right)(1-q)^{k-j}(1-p)^{k-j}(-r)^j.
$$
This completes the proof.

\section{Concluding remarks}

Formulae for the distribution of occupancy times obtained in Section~1
are probably new, but of limited practical use, as they require explicit
information about powers of $P_{U^cU^e}$; in the two-state case all the
relevant ingredients are scalars. In addition to the work of 
Sericola~\cite{Ser00}, already mentioned, there is also a 
connection with first-passages
times, which is well known and exploited in the classification of
states (see, for example, Theorem~2.5.9 of~\cite{Bre20}). On the other
hand, {\em expected\/} occupancy times have been studied widely.
Especially important is the connection with potential theory (see
Section~4.2 of~\cite{Nor97}) where, in the present discrete-time
setting, (\ref{Nn}) would be replaced by $F_n = \sum_{m=1}^n f(X_m)$,
where $f(i)$ is a cost associated with being in state~$i$, and $F_n$ is
total cost up to time $n$.
It would be interesting to see if any of the results presented here
could be extended to this more general setting.
However, even expected occupancy times cannot 
be obtained readily from~(\ref{Gtk}). An alternative approach is needed:
consider the probability generating function of $N_n$, conditional on $X_0=i$,
$$
H_i(n,z)=\E (z^{N_n}|X_0=i)=\sum_{k=0}^n g_i(n,k) z^k,
\qquad |z|<1,
$$
and the corresponding column vector $H(z)= (H_i(z),\, i\in S)$ given by
$H(n,z)= \sum_{k=0}^n g(n,k) z^k$.
Then, (\ref{gnk}) and (\ref{gn0}) give
\begin{equation}
H(n,z)=(B+Az)H(n-1,z), \qquad n\geq 2, 
\label{Hnz}
\end{equation}
with $H(1,z)=(B+Az)\bone$, 
and hence $H(n,z)=(B+Az)^n\bone$, $n\geq 0$. However, since $A$ and $B$
do not commute, taking derivatives near $z=1$ is
not straightforward, and it is better to deal with~(\ref{Hnz}) when
evaluating moments of $N_n$. Differentiating~(\ref{Hnz}), and
letting $z\uparrow 1$, gives us an expression for 
$e(n)=(e_i(n),\, i\in S)$, where $e_i(n)=\E(N_n|X_0=i)$.
We find that $e(n)= Pe(n-1) + A\bone$, $n\geq 2$, with
$e(1)=A\bone$, and hence that
$$
e(n)=(I+P+\dots +P^{n-1})A\bone,
\qquad n\geq 1.
$$
We might also be able to prove results on the total cost 
up to some terminal time (for instance the time of extinction),
results analogous to those presented in \cite{Pol03,PS02a} for Markov chains
in continuous-time.

\bibliographystyle{plain}

\def\cprime{$'$}

\end{document}